\documentclass[11pt]{article}
\usepackage{epsfig}

\usepackage{amsfonts,amssymb,theorem,amsmath,amscd}
\usepackage[mathscr]{eucal}
      \makeatletter
\def\setseccntfmt{\renewcommand{\@seccntformat}[1]{\S
    \csname the##1\endcsname.\hspace{1ex}}}
\def\setsubseccntfmt{\renewcommand{\@seccntformat}[1]{%
    (\csname the##1\endcsname)\hspace{0.5ex}}}
\renewcommand{\section}{\setseccntfmt\@startsection
  {section}{1}{0mm}{-\baselineskip}{0.5\baselineskip}{\sf\bfseries\large}}
\def\presubsection{\setsubseccntfmt\@startsection
  {subsection}{2}{0mm}{-\baselineskip}{-0.5ex}{\bfseries\upshape}}
\def\tmpa{}\def\tmpb{}
\newcommand{\addspaceifnonempty}[1]{\def\tmpa{}\def\tmpb{#1}%
  \ifx\tmpa\tmpb{}\else{\hspace{0.5ex}#1\hspace{1ex}}\fi}
\renewcommand{\subsection}[1][]{\presubsection{\addspaceifnonempty{#1}}}
\def\presubsubsection{\setsubseccntfmt\@startsection
  {subsubsection}{3}{0mm}{-\baselineskip}{-0.5ex}{\bfseries\upshape}}
\renewcommand{\subsubsection}[1][]{\presubsubsection{\addspaceifnonempty{#1}}}
\gdef\th@nonumplain{\normalfont\itshape
  \def\@begintheorem##1##2{\item[\hskip\labelsep\theorem@headerfont ##1]}%
  \def\@opargbegintheorem##1##2##3{%
    \item[\hskip\labelsep \theorem@headerfont ##1\ (##3)]}}
\gdef\th@change{
  \def\@begintheorem##1##2{\item[\hskip\labelsep
    {\bfseries\upshape(##2)\hskip 1ex}\theorem@headerfont ##1]}%
  \def\@opargbegintheorem##1##2##3{\item[\hskip\labelsep
    {\bfseries\upshape(##2)\hskip 1ex}\theorem@headerfont ##1\ (##3)]}}
\makeatother   \def\nthm{\newtheorem}
\theoremheaderfont{\normalfont\bfseries\scshape}
{\theoremstyle{change}{\theorembodyfont{\normalfont\itshape}
\nthm{conj}[subsection]{Conjecture}  \nthm{*conj}[subsubsection]{Conjecture}
\nthm{thm}[subsection]{Theorem}      \nthm{*thm}[subsubsection]{Theorem}
\nthm{prop}[subsection]{Proposition} \nthm{*prop}[subsubsection]{Proposition}
\nthm{lemma}[subsection]{Lemma}      \nthm{*lemma}[subsubsection]{Lemma}
\nthm{cor}[subsection]{Corollary}    \nthm{*cor}[subsubsection]{Corollary}}
{\theorembodyfont{\normalfont\rmfamily}
\nthm{rem}[subsection]{Remark}       \nthm{*rem}[subsubsection]{Remark}
\nthm{defn}[subsection]{Definition}  \nthm{*defn}[subsubsection]{Definition}}}
{\theoremstyle{nonumplain}
{\theorembodyfont{\normalfont\itshape}
\nthm{thm*}[subparagraph]{Theorem}   \nthm{prop*}[subparagraph]{Proposition}
\nthm{lemma*}[subparagraph]{Lemma}   \nthm{cor*}[subparagraph]{Corollary}
\nthm{conj*}[subparagraph]{Conjecture}}
{\theorembodyfont{\normalfont\rmfamily}
 \nthm{rem*}[subparagraph]{Remark}   \nthm{defn*}[subparagraph]{Definition}}}
\DeclareMathAlphabet\eusm{U}{eus}{m}{n}

\def\makeop#1{\expandafter\def\csname#1\endcsname
  {\mathop{\rm #1}\nolimits}\ignorespaces}
\makeop{Hom}   \makeop{End}   \makeop{Aut}   \makeop{Isom}  \makeop{Pic} 
\makeop{Gal}   \makeop{ord}   \makeop{Char}  \makeop{Div}   \makeop{Lie} 
\makeop{PGL}   \makeop{Corr}  \makeop{PSL}   \makeop{sgn}   \makeop{Spf}
\makeop{Spec}  \makeop{Tr}    \makeop{Nr}    \makeop{Fr}    \makeop{disc}
\makeop{Proj}  \makeop{supp}  \makeop{ker}   \makeop{im}    \makeop{dom}
\makeop{coker} \makeop{Stab}  \makeop{SO}    \makeop{SL}    \makeop{SL}
\makeop{Cl}    \makeop{cond}  \makeop{Br}    \makeop{inv}   \makeop{rank}
\makeop{id}    \makeop{Fil}   \makeop{Frac}  \makeop{GL}    \makeop{SU}
\makeop{Trd}   \makeop{Sp}    \makeop{Tr}    \makeop{Trd}   \makeop{diag}
\makeop{Res}   \makeop{ind}   \makeop{depth} \makeop{Tr}    \makeop{st}
\makeop{Ad}    \makeop{Int}   \makeop{tr}    \makeop{Sym}   \makeop{can}
\makeop{length}\makeop{SO}    \makeop{torsion}
\def\makebb#1{\expandafter\def
  \csname bb#1\endcsname{{\mathbb{#1}}}\ignorespaces}
\def\makebf#1{\expandafter\def\csname bf#1\endcsname{{\bf #1}}\ignorespaces}
\def\makegr#1{\expandafter\def
  \csname gr#1\endcsname{{\mathfrak{#1}}}\ignorespaces}
\def\makescr#1{\expandafter\def
  \csname scr#1\endcsname{{\mathscr{#1}}}\ignorespaces}
\def\makecal#1{\expandafter\def\csname cal#1\endcsname{{\cal #1}}\ignorespaces}
\def\doLetters#1{#1A #1B #1C #1D #1E #1F #1G #1H #1I #1J #1K #1L #1M
                 #1N #1O #1P #1Q #1R #1S #1T #1U #1V #1W #1X #1Y #1Z}
\def\doletters#1{#1a #1b #1c #1d #1e #1f #1g #1h #1i #1j #1k #1l #1m
                 #1n #1o #1p #1q #1r #1s #1t #1u #1v #1w #1x #1y #1z}
\doLetters\makebb   \doLetters\makecal  \doLetters\makebf \doLetters\makescr
\doletters\makebf   \doLetters\makegr   \doletters\makegr
     \def\qed{\qedmark\medbreak}%
\def\qedmark{{\enspace\vrule height 6pt width 5pt depth 1.5pt}}%

\newcommand{\rmitem}[1][{}]{\item[{\rm#1}]}

\makeatletter   \newdimen\mina@@\mina@@=18pt
\newcommand{\xrtarw}[2][]{\mathrel{\mathop{\,\setbox\z@\vbox{\m@th
  \hbox{$\scriptstyle\;{#1}\;\;$}\hbox{$\m@th\scriptstyle\;{#2}\;\;$}}%
  \hbox to\ifdim\wd\z@>\mina@@\wd\z@\else\mina@@\fi{\rightarrowfill@
  \displaystyle}\,}\limits^{#2}\@ifnotempty{#1}{_{#1}}}}
\newcommand{\xltarw}[2][]{\mathrel{\mathop{\,\setbox\z@\vbox{\m@th
  \hbox{$\scriptstyle\;\;{#1}\;$}\hbox{$\m@th\scriptstyle\;\;{#2}\;\;$}}%
  \hbox to\ifdim\wd\z@>\mina@@\wd\z@\else\mina@@\fi{\leftarrowfill@
  \displaystyle}\,}\limits^{#2}\@ifnotempty{#1}{_{#1}}}}
\makeatother
\def\XYmatrix{\xymatrix@M=8pt} 
\def\ncmd{\newcommand}
\ncmd{\xysubset}[1][r]{\ar@<-2.5pt>@{^(-}[#1]\ar@<2.5pt>@{_(-}[#1]}
\ncmd{\XYmatrixc}[1]{\vcenter{\XYmatrix{#1}}}
\ncmd{\xyto}[1][r]{\ar@{->}[#1]}      \ncmd{\xyinj}[1][r]{\ar@{^(->}[#1]}
\ncmd{\xysurj}[1][r]{\ar@{->>}[#1]}   \ncmd{\xyline}[1][r]{\ar@{-}[#1]}
\ncmd{\xydotsto}[1][r]{\ar@{.>}[#1]}  \ncmd{\xydots}[1][r]{\ar@{.}[#1]}
\ncmd{\xyleadsto}[1][r]{\ar@{~>}[#1]} \ncmd{\xyeq}[1][r]{\ar@{=}[#1]}
\ncmd{\xyequal}[1][r]{\ar@{=}[#1]}    \ncmd{\xyequals}[1][r]{\ar@{=}[#1]}
\ncmd{\xymapsto}[1][r]{l\ar@{|->}[#1]}\ncmd{\xyimplies}[1][r]{\ar@{=>}[#1]}

\normalsize

\makeop{Bl}

\providecommand{\rmitem}[1][{}]{\item[\rm#1]}

\def\Spec{{\rm Spec}}

\def\Qbar{\overline{\bbQ}}

\newcommand{\Z}{\mathbb Z}
\newcommand{\Q}{\mathbb Q}
\newcommand{\R}{\mathbb R}
\newcommand{\C}{\mathbb C}
\renewcommand{\H}{\mathbb H}  
\renewcommand{\O}{\mathcal O} 
\newcommand{\F}{\mathbb F}

\newcommand{\npr}{\noindent }

\newcommand{\<}{\langle}   
\renewcommand{\>}{\rangle} 

\newcommand{\ch}{characteristic }
\newcommand{\ac}{algebraically closed }
\newcommand{\dieu}{Dieudonn\'{e} }

\ncmd{\ol}{\overline}
\ncmd{\wt}{\widetilde}
\ncmd{\isoto}{\stackrel{\sim}{\to}}
\ncmd{\ul}{\underline}

\makeop{Ram}
\makeop{Rep}
\makeop{Norm}

\begin{document}
\begin{center}
\Large\scshape
On the slope stratification of certain Shimura
varieties
\end{center}
\medbreak

\centerline{\scshape Chia-Fu Yu\footnote{\today} 
%
}

\begin{abstract}
  In this paper we study the slope stratification on the good
  reduction of the type C family Shimura varieties. 
  We show that there is an open dense subset $U$ of the
  moduli space such that any point in $U$ can be deformed to a point
  with a given lower {\it admissible} Newton polygon. For the Siegel
  moduli spaces, this is obtained by F.~Oort which plays an important
  role in his proof of the strong Grothendieck conjecture concerning
  the slope stratification. 
  We also investigate the $p$-divisible groups and their isogeny classes
  arising from the abelian varieties in question.
\end{abstract}
\section{Introduction}
\label{sec:01}

Let $p$ be a rational prime number. Let $B$ be a finite-dimensional division
algebra over $\Q$ with a positive involution $*$. We consider the
cases where $B$ is either 

(a) a totally real field $F$ of degree $d$, or

(b) a totally indefinite quaternion algebra over a totally real field
$F$ of degree $d$.

Let $O_B$ be an order of $B$ stable under the involution such that
$O_B\otimes \Z_p$ is a maximal order of $B\otimes \Q_p$.
Let $\calM$ denote the moduli stack over $\ol \F_p$ of separably
polarized abelian $O_B$-varieties of dimension $g$, where $g=dn$ in
case (a) and $g=2dn$ in case (b) for some positive integer $n$. 
A polarized abelian $O_B$-variety $(A,\lambda,\iota)$ is a polarized
abelian variety $(A,\lambda)$ together with a ring 
embedding $\iota:O_B\to \End(A)$ such that 
$\lambda \iota(b^*)=\iota(b)^t \lambda$ for all
$b\in O_B$. For a polarized abelian $O_B$-variety $\ul
A=(A,\lambda,\iota)$, the associated $p$-divisible group $\ul H=\ul
A[p^\infty]$ with the additional structure is a quasi-polarized
$p$-divisible $O_B\otimes \Z_p$-group (in \cite{yu:reduction}), or a
quasi-polarized BT $O_B\otimes \Z_p$-module (in \cite{moonen:eo}). 
Define the (refined) Newton polygon of $\ul H$ by 
\[ NP(\ul H)=(NP(H_v))_{v|p}, \]
where $H=\oplus_v H_v$ is the decomposition of $p$-divisible groups 
coming from $O_F\otimes \Z_p=\oplus_v O_v$ 
and $NP(H_v)$ is the usual Newton polygon of $H_v$. A natural question
is whether a given quasi-polarized $p$-divisible $O_B\otimes
\Z_p$-group over an \ac field $k$ comes from a quasi-polarized abelian
$O_B$-variety. Let $M$ be the covariant \dieu module associated to
$\ul H$ and let $M=\oplus_{v|p}M_v$ the corresponding decomposition. 
One can show that each factor $M_v$, considered as an $O_v\otimes
W(k)$-module, is 
free. If $M$ arises from an abelian $O_B$-variety, then it satisfies
an additional condition that the ranks $\rank_{O_v \otimes W(k)} M_v$
are the same for all $v|p$. The latter is equivalent to the following 
condition:
\[ (*)\quad \Tr_{W(k)} (a;M)\cdot [F:\Q]=(\rank_{W(k)} M)\cdot
\Tr_{F/\Q}(a), \quad 
\text{for all $a\in O_F$}. \] 
In this paper we prove

\begin{thm}\label{11}
Let $k$ be an \ac field of \ch $p$.
\rmitem[(1)] Any quasi-polarized supersingular $p$-divisible $O_B\otimes
\Z_p$-group over $k$ whose \dieu module satisfies ($*$) is realized by
a quasi-polarized abelian 
$O_B$-variety $\ul A$.
\rmitem[(2)] Let $\ul H$ be a quasi-polarized $p$-divisible $O_B\otimes
\Z_p$-group over $k$. If its Newton polygon $NP(\ul H)$ is realized by
an object $\ul A'$, then $\ul H$ is realized by an object $\ul A$. 
\end{thm}

The statement (2) is an extension of (a special case of) a result of
Rapoport and Richartz \cite[Proposition 1.17]{rapoport-richartz} on
the Newton map. The statement (1) follows from (2) and a construction
of supersingular points; see Theorem~\ref{25} for the statement. 
The existence of supersingular points is well-known, at least for some 
special cases. We simply supply a more complete reference here. \\

A Newton polygon $\ul \beta=(\beta_v)_v$ is
called {\it admissible} if it arises from a quasi-polarized
$p$-divisible $O_B\otimes \Z_p$-group $\ul H$ satisfying ($*$). 
The partial order imposed on the set of admissible Newton polygons
here is given by 
$(\beta_v)_v\ge (\beta'_v)_v$ if every component $\beta_v$ lies
above or equal to $\beta'_v$. In other words, the associated Newton
stratum with higher Newton polygon has smaller dimension. 
For the readers who are familiar with Shimura varieties, the admissible
Newton polygons correspond exactly the elements in the image of the
set $B(G,\mu)$ (e.g. in \cite[Section 4]{rapoport:modp}) of
$\mu$-admissible $\sigma$-conjugacy classes under the 
Newton map, where $G$ is the base change to $\Q_p$ of the group defining
Shimura varieties and $\mu$ is the corresponding dominant minuscule
coweight. 

Let $\Delta$ be the product of rational primes at which
$B$ is ramified. 
   
\begin{thm}\label{12}
  Assume that $p \nmid \Delta$. Let $U$ be an open subset of $\calM$ defined
  by 
\[ U=\left\{ \ul A=(A,\lambda,\iota)\in \calM\, |\, a(H_v)\le 1, \forall\,
  v|p\,\right\}, \]
  where $H_v$ is the $v$-component of $\ul A[p^\infty]$, $a(H_v):=\dim_k
\Hom(\alpha_p, H_v)$ is the $a$-number of $H_v$, and $v$ runs through
the primes of $O_F$ over $p$. 

Then for any point $u\in U$ with Newton polygon $\ul \beta$ and for
any admissible Newton polygon $\ul \beta'$ below $\ul \beta$ there is
a generalization $u'\in \calM$ of $u$ whose Newton polygon is equal to
$\ul \beta'$. 
\end{thm}

The open subset $U$ is dense in $\calM$ as it contains the ordinary locus.
When $B=\Q$, Theorem~\ref{12} is proved by F.~Oort 
\cite{oort:cayley}, which plays an important role in his proof of the 
{\it strong Grothendieck conjecture} for the Siegel moduli space
\cite{oort:grothendieck}: replacing $U$ by the total space $\calM$ in
Theorem~\ref{12}. See Rapoport \cite[Section 5]{rapoport:np} and 
\cite[Section 7]{rapoport:modp} for further discussions.  
For the present cases, this result is not enough to 
conclude the analogue of the strong Grothendieck conjecture, as the
intersection of a Newton stratum $\calW$ with $U$ is not dense in
$\calW$ in general \cite{yu:thesis}. 
Nevertheless, Theorem~\ref{12} concludes a few new results 
and improves some previous results in the type C family of 
Shimura varieties.  

\begin{thm}\label{13}
  Assume that $p \nmid \Delta$. For any sequence of admissible Newton
  polygons $\ul {\beta_1} > \ul {\beta_2} > \dots > \ul {\beta_s}$, there
  is a chain of irreducible substacks $V_1\subset V_2\subset\dots
  \subset V_s$ of $\calM$ such that the Newton polygon $NP(\ul
  A_{\eta_i})$ of the generic point $\eta_i$ is $\ul {\beta_i}$. 
\end{thm}

This is a weak form of the analogue of the strong Grothendieck
conjecture (or called the {\it weak Grothendieck conjecture} in
Oort \cite{oort:cayley}). An immediate consequence of this is that
every Newton stratum is non-empty. This proves a special case
of Conjecture~7.1 in Rapoport \cite{rapoport:modp}. 

\begin{thm}\label{14}
  When $p \nmid \Delta$, any quasi-polarized $p$-divisible
  $O_B\otimes \Z_p$-group over an \ac field of \ch $p$ satisfying
  ($*$) is realized by a quasi-polarized abelian $O_B$-variety.
\end{thm}

When $\calM$ is smooth, Theorem~\ref{14} shows that any {\it $p$-adic
invariant} stratum arising from $p$-divisible groups in question is
non-empty. This improves previous results on stratifications in a
special case (type C) by removing the 
non-emptiness assumption; for example, see the main theorem (Corollary
3.1.6) in Moonen \cite{moonen:eo}.  
One would expect that Theorem~\ref{14} holds as well without the
condition $p \nmid \Delta$. This is a technical point as the proof
presented here replies on an explicit deformation computation, which
is more involved in the ramified cases. We hope to return to this
technical issue soon.   

Finally we remark that the strong Grothendieck conjecture has been
proved by B\"ultel and Wedhorn \cite{bueltel-wedhorn} in the case of
inner forms of groups ${\rm U}(n,1)$ associated to imaginary quadratic
fields at unramified primes $p\not =2$, where the split case has been
known for a long time as the Lubin-Tate formal moduli. 

The paper is organized as follows. In Section 2 we construct a
supersingular abelian variety with the additional structure; the basic
tool is the classification of positive involutions up to conjugate. In
Section 3 we show that the formal isogenies in question are determined by
the Newton polygons. In Sections 2 and 3 there is no assumption of
$p$. For the rest (Sections 4 and 5) we assume that $p\nmid \Delta$.
Using the Morita equivalence, one reduces to considering the deformation
problem on local components of $p$-divisible groups in question. We
treat the reduced local situation directly. In Section 4 we classify
the \dieu modules, with the additional structure, with $a$-number one. 
In Section 5 we recall the Cartier-\dieu theory and the method of
Cayley-Hamilton. Then we construct the universal deformation and compute
any formal completion of any Newton stratum on the open subset $U$ in
Theorem~\ref{12}.

An experienced reader would find that the present formulation differs
a bit from the literatures \cite{kottwitz:isocrystals,
rapoport-richartz, chai:np, rapoport:modp}. The choice is simply
because it is more direct to adopt the explicit deformation theory of
displays and the Cayley-Hamilton method. As the groups defining
Shimura varieties treated here satisfies the Hasse principle, the
translation is straightforward. \\

\npr {\it Acknowledgments.} The author is deeply indebted to Frans Oort
  for his inspiring work \cite{oort:cayley}. An earlier draft (MPI-preprint)
  of this paper was done during the author's stay at the
  Max-Planck-Institut f\"ur Mathematik in Bonn in the 2001 Summer. He
  wishes to thank the Institute for the kind hospitality and the
  excellent working environment. One should mention that Rapoport's
  article \cite{rapoport:modp} has been very helpful when the author 
  reorganizes the present paper. Finally, he acknowledges the
  referee for his encourgements and helpful comments.   


\section{Construction of supersingular points}
\label{sec:02}

In the rest of this paper, $k$ denotes an \ac field of \ch $p$. 
We keep the notation of the previous section. 
In this section we will construct a polarized supersingular
abelian $O_B$-variety of dimension $g$ over $k$, where $g=dn$ in case (a) and
$g=2dn$ in case (b) for some positive integer $n$.   

\subsection{}
\label{21}
We recall basic results on simple algebras over $\Q$ with
involutions (cf. \cite{kottwitz:jams92}). Let $D$ be a
finite-dimensional simple algebra with a positive involution $*$. 
Then we have $(D\otimes_\Q \Qbar, *)\simeq \oplus_{i=1}^r (L, *_1)$, 
where $r=[{\rm Cent}(D)^{*=1}:\Q]$ and $(L,*_1)$ is one of the
following three cases: 
\begin{itemize}
\item [(A)]: $L=M_n(\Qbar)\oplus M_n(\Qbar)$ and $(A,B)^{*_1}=(B^t,A^t)$,
\item [(C)]: $L=M_{2n}(\Qbar)$ and $*_1$ is the standard symplectic
  involution,
\item [(BD)]: $L=M_n(\Qbar)$ and $*_1$ is the transposition.
\end{itemize}
We also have the corresponding cases for $(D\otimes_\Q \R, *)$:
$(D\otimes_\Q \R, *)\simeq \oplus_{i=1}^r(N,*_2)$ for the same $r$, 
where $(N,*_2)$ is one of the following
\begin{itemize}
\item [(A)]: $N=M_n(\C)$ and $(a_{ij})^{*_2}=(\bar{a}_{ji})$,
\item [(C)]: $N=M_n(\H)$ and $(a_{ij})^{*_2}=(\bar{a}_{ji})$,
\item [(BD)]: $N=M_n(\R)$ and $(a_{ij})^{*_2}=({a}_{ji})$.
\end{itemize}
Here $a\mapsto \bar a$ is the standard conjugation.

\subsection{}
\label{22}
Let $\,'$ be another positive involution on $D$.
Then there is an positive element $c\in D$ with $c^*=c$ such that
$x'=cx^* c^{-1}$ for all $x\in D^\times$ \cite[Lemma
2.11]{kottwitz:jams92}. Let $U(D,*)_\Q$ denote 
the associated unitary group, which associates to a commutative
$\Q$-algebra $R$ the group of $R$-points
\[ U(D,*)(R):=\{x\in (D\otimes_\Q R); x x^*=1 \}. \]
Define the algebraic variety $X_c$ over $\Q$ by
\[ X_c(R):=\{x\in (D\otimes_\Q R); x x^*=c \}. \]
Note that $X_c(R)\not=\emptyset$ if and only if $(D\otimes R, ')\simeq
(D\otimes R, *)$. One has $X_c(\Qbar)\not =\emptyset$ and the group 
$U(D,*)$ acts simply-transitively on $X_c$ from the right. 
In other words, $X_c$ is a right principal homogeneous
space under the group $U(D,*)$. Therefore it gives rise to an 
element $[c]$ in ${\rm H}^1(\Q, U(D,*))$. Note that the image of 
$[c]$ in ${\rm H}^1(\R, U(D,*))$ is trivial. 

\begin{lemma}\label{23} Notations as above. If $(D,*)$ is of type (C),
  then $(D,\,')\simeq (D,*)$ over $\Q$. 
\end{lemma}
\begin{proof}
  In this case, the group $U(D,*)$ is semi-simple and 
  simply-connected. It follows from Kneser's theorem and the Hasse
  principle that $H^1(\Q,U(D,*))=0$, hence $(D,\,')\simeq (D,*)$ over
  $\Q$. \qed
\end{proof}

\subsection{}
\label{24}
Let $(D_1,*_1)$ and $(D_2,*_2)$ be finite-dimensional simple algebras
over $\Q$ with positive involutions. Suppose that the fixed fields of
the centers are the same, called $F$. For simplicity, assume that 
${\rm Cent}(D_1)\otimes_F {\rm Cent}(D_2)$ is a field. We define a new 
simple algebra with a positive involution by 
\[ (D,*):=(D_1\otimes_F D_2, *), \]
where $(d_1\otimes d_2)^*:=d_1^{*_1}\otimes d_2^{*_2}$. 
Clearly $*$ is an involution on $D$. To see its positivity, 
first we have the following formula, for $a\in D_1$ and $b\in D_2$,
\begin{equation*}
  \label{eq:241}
  \begin{split}
    \Tr _{D_1\otimes_F D_2/\Q}(a\otimes b)& =\sum_{\sigma:F\hookrightarrow
  \Qbar} \Tr_{[(D_1\otimes_{F,\sigma}\Qbar)\otimes
  (D_2\otimes_{F,\sigma} \Qbar)] /\Qbar}
  (\sigma(a)\otimes \sigma(b))\\
  &= \sum_{\sigma:F\hookrightarrow \Qbar}
  \Tr_{(D_1\otimes_{F,\sigma}\Qbar)/\Qbar}(\sigma(a))\cdot 
  \Tr_{(D_2\otimes_{F,\sigma} \Qbar) /\Qbar}(\sigma(b)) 
  \end{split}
\end{equation*}
It follows that for $a_1,a_2\in D_1$ and $b_1,b_2\in D_2$,
\begin{equation*}
  \begin{split}
    \Tr _{D_1\otimes_F D_2/\Q}(&(a_1\otimes b_1)(a_2\otimes b_2)^*)\\
   & =
   \sum_{\sigma:F\hookrightarrow \R}
  \Tr_{(D_1\otimes_{F,\sigma}\R)/\R}(\sigma(a_1)\sigma(a_2)^{*_1})\cdot 
  \Tr_{(D_2\otimes_{F,\sigma} \R) /\R}(\sigma(b_1)\sigma(b_2)^{*_2}), 
  \end{split}
\end{equation*}
which is a positive form. 
Moreover, we have the following results:
\begin{itemize}
\item If one of them is of type (A), then $(D,*)$ is of type (A).
\item If both are of type (C) or of type (BD), then $(D,*)$ is of type (BD).
\item If one is of type (C) and the other is of type (BD), then $(D,*)$ is
  of type (C).
\end{itemize}

\begin{thm}\label{25} 
  There exists a polarized supersingular abelian $O_B$-variety of
  dimension $g$ over $k$. 
\end{thm}
\begin{proof} Case (a): $B=F$, $[F:\Q]=d$ and $g=dn$. 
  Choose a supersingular elliptic curve $E$ and put $A_0:=(E\otimes_\Z
  O_F)^{n}$. One can choose a polarization on $A_0$, induced from
  one on $E$, and make $A_0$ a polarized supersingular abelian
  $O_F$-variety.    

  Case (b): $[B:F]=4$, $[F:\Q]=d$ and $g=2dn$. 
  Choose a supersingular elliptic curve $E$ and put $A_0:=(E\otimes
  O_F)^{2n}$ as before. Let $(D,\,')$ be 
  the simple algebra $\End_{O_F}(A_0)\otimes \Q$ with 
  a Rosati involution. Note that $D\simeq M_{2n}(D_F)$, where
  $D_F:=D_p\otimes_\Q F$ and $D_p$ is the quaternion algebra
  over $\Q$ ramified exactly at $\{\infty, p\}$. The nontrivial 
  Brauer invariants for $D_F$
  are $\frac{1}{2}$ at each infinity and $\frac{1}{2}$ at some
  primes above $p$. Let $C$ be the quaternion algebra over $F$ such
  that the Brauer invariants of $B\otimes_F C$ are the same as
  those of $D$. Replacing $C$ by $M_n(C)$, we have $B\otimes_F
  C\simeq D$. Choose any positive involution on $C$ and make the
  product positive involution $*$ on $B\otimes_F C$, which extends 
  $*$ on $B$. By Lemma~\ref{23}, there
  is an isomorphism $(B\otimes_F C,*)\simeq (D,\,')$. Therefore, we 
  construct a polarized supersingular abelian $O_F$-variety $A_0$
  together with a monomorphism $B\to
  \End_{O_F}(A_0)\otimes \Q$ preserving the involutions. 
  Finally, one can choose an abelian
  $O_B$-variety in the $B$-linear isogeny class of $A_0$ and then 
  take the induced $B$-linear polarization. 
  The construction is complete. \qed
\end{proof}

\begin{rem}
  There is another way to construct supersingular points as in
  (\ref{25}). One first constructs a supersingular abelian variety up
  to isogeny together with a $B$-action. This can be done by choosing
  an embedding from $B$ to the endomorphism algebra. 
  Then a result of Kottwitz \cite[Lemma 9.2]{kottwitz:jams92}
  furnishes a $B$-linear polarization. This construction has been used
  in B\"ultel and Wedhorn \cite[Lemma 5.2]{bueltel-wedhorn}. 
\end{rem}

\section{Formal isogenies}
\label{sec:03}

We keep the notation as in the previous sections. Theorem~\ref{11} (2)
follows from the following

\begin{prop}\label{31}
  Let $\ul A_1, \ul A_2$ be polarized abelian $O_B$-varieties over
  $k$. If $NP(\ul A_1)=NP(\ul A_2)$, then there is a polarized
  $O_B\otimes \Z_p$-linear quasi-isogeny from $\ul A_1[p^\infty]$ to $\ul
  A_2[p^\infty]$ 
\end{prop}

This means that the quasi-isogeny $\varphi\in
\Hom(A_1[p^\infty],A_2[p^\infty])\otimes \Q$ is
$O_B\otimes\Z_p$-linear and that $\varphi^*\lambda_2=\lambda_1$.
Since the situation is local, one will replace $F$ by its completion
at one prime and $B$ by its completion. Let $H_1$ and $H_2$ be the
corresponding factor of $\ul A_1[p^\infty]$ and $\ul A_2[p^\infty]$
respectively. We have two cases, by the Morita equivalence:

(a) $O_B=O_F$. 

(b) $B$ is a division quaternion algebra over $F$ with an involution
$x^*=\alpha x' \alpha^{-1}$ for some $\alpha\neq 0$ and 
$\alpha'=-\alpha$, where $x\mapsto x'$ is the standard involution.  

\begin{lemma}\label{32}
  Let $D$ be a finite-dimensional $\Q_p$-simple algebra, and $O_D$ a
  maximal order of $D$. Then any two $p$-divisible $O_D$-groups
  $H_1,H_2$ with $NP(H_1)=NP(H_2)$ are $O_D$-linearly isogenous.
\end{lemma}
\begin{proof}
The proof is elementary; one can use the Noether-Skolem theorem to
adjust a chosen isogeny to be $O_D$-linear (cf. \cite[Section
10]{yu:thesis}), or use a result of Rapoport and Zink (see
Remark~\ref{34} below) to conclude this. Hence we omit the proof. \qed   


\end{proof}

Now we choose an $O_B$-linear isogeny $\varphi:H_1\to H_2$. 
The pullback of $\lambda_2$ has the form $\lambda_1 a$ for some $a\in
\End_{O_B}(H_1)\otimes \Q_p$ which is invariant under the Rosati
involution $*_1$ by $\lambda_1$. The element $[a]$ defines
$J_1$-torsor over $\Q_p$, where $J_1$ is the unitary algebraic group 
$\Aut_{B}(H_1,*_1)$ over $\Q_p$. This $J_1$-torsor vanishes if and only
if there is an element $b\in \End_{O_B}(H_1)\otimes \Q_p$ such that
$b^{*_1}b=a$. Note that for $b\in \End_{O_B}(H_1)$, the pull-back
$b^*\lambda_1$ is $\lambda_1 b^{*_1} b$.
By Lemma~\ref{33} below, there is an
element $b\in \End_{O_B}(H_1)\otimes \Q_p$ such that $b^{*_1}
b=a^{-1}$. Then one shows $(\varphi b)^*\lambda_2=\lambda_1$. This
finishes the proof of Proposition~\ref{31}. 

\begin{lemma}\label{33}
  Let $J_1$ be as above. One has ${\rm H}^1(\Q_p,J_1)=0$. 
\end{lemma}
\begin{proof}
Let $[H]$ denote the isogeny class of a $p$-divisible group $H$.
We write 
\[ [H_1]=\oplus_{\lambda<1/2}[H_{1,\lambda}\oplus
H_{1,1-\lambda}]\oplus [H_{1,1/2}] \]
into the slope decomposition. The quasi-polarization sends each factor
onto its dual.  One has 
\[ J_1=\prod_{\lambda<1/2} J_{1,\lambda}\, \times\,
  J_{1,1/2}, \]
where $J_{1,\lambda}$ (resp. $J_{1,1/2}$) is the unitary group
$\Aut_{B}([H_{1,\lambda}\oplus H_{1,1-\lambda}],*_1)$ (resp. $\Aut
 _{B}([H_{1,1/2}],*_1)$) over $\Q_p$. Therefore,
${\rm H}^1(\Q_p,J_1)=\prod_{\lambda<1/2} {\rm H}^1
(\Q_p,J_{1,\lambda})\, \times {\rm H}^1(\Q_p,J_{1,1/2})$. 

The group $J_{1,\lambda}$ is the algebraic group associated to the
multiplicative group $\End_B([H_{1,\lambda}])^\times$. By Hilbert's 90
theorem, ${\rm H}^1(\Q_p,J_{1,\lambda})=0$. As the group $J_{1,1/2}\otimes \ol
\Q_p$ is isomorphic to a product of symplectic groups, 
$J_{1,1/2}$ is semi-simple and simply connected. This shows
${\rm H}^1(\Q_p,J_{1,1/2})=0$ and completes the proof. \qed   
\end{proof}

\begin{rem}\label{34}
   In the proof (\ref{31}) we proved a fact that the fibers of the Newton
   map are classified by the Galois cohomology ${\rm H}^1(\Q_p,J)$ for
   the twisted centralizer $J$. This result has 
   been known due to Rapoport and Richartz \cite{rapoport-richartz}
   for general F-isocrystals with additional structures. Here we simply
   give an elementary presentation of this result for our special
   cases. Rapoport and Zink show that $J$ is an inner form of the
   centralizer of $s\nu(p)$ in the structure group $G_{\Q_{p^s}}$, 
   where $\nu$ is the attached slope homomorphism (see
   \cite[1.7--1.16]{rapoport-zink} for details). From this one can
   determine $J$ and compute ${\rm H}^1(\Q_p, J)$ as well. 
\end{rem}

\section{The main theorem}
\label{sec:04}

\def\iz{i\in I}
\subsection{} \label{41}
In this and the next sections, $F$ denotes, following the convention,
the Frobenius operator on a \dieu module. Let $W:=W(k)$ be the 
ring of Witt vectors, and $\sigma$ the Frobenius map on $W$. 
Let $\O$ be an \'etale
extension of $\Z_p$ of degree $f$, that is, $\O\simeq W(\F_{p^f})$. 
A quasi-polarized \dieu $\O$-module is a \dieu
module $M$ together with a non-degenerate alternating $W$-linear pairing
(quasi-polarization) $\<\, ,\, \>: M\times M\to W$ such that
$\<Fx,y\>=\<x,Vy\>^\sigma$ for all $x,y\in M$, and with a ring embedding
$\iota: \O\to \End_{\rm dieu}(M)$ such that $\<ax,y\>=\<x,ay\>$ for
all $a\in \O$, $x,y \in M$. 

Let $I:=\Hom(\O,W)$ be the set of embeddings. We may and do choose
an identification $\Z/f\Z\simeq I$, mapping $i\mapsto \sigma_i$,
so that $\sigma \sigma_i=\sigma_{i+1}$.

Let $\ul M=(M,\<\, ,\, \>,\iota)$ be a quasi-polarized \dieu
$\O$-module over $k$. 
Denote by
\[ M^i:=\{x\in M\, |\, a x=\sigma_i(a) x, \ \forall\ a\in
\O \, \} \]  
the $\sigma_i$-component of $M$\!, which is a free $W$-module of
rank $2r$. We have the decomposition 
\[ M=M^0\oplus M^1\oplus \dots \oplus M^{f-1} \]
in which  $F:M^i\to M^{i+1}$, $V:M^{i+1}\to M^i$. 
The summands $M^i, M^j$ are orthogonal with respect to the pairing
$\<\, ,\,\>$ for $i\not=j$. Conversely, 
a \dieu module together with such a decomposition and these properties
is a quasi-polarized \dieu $\O$-module. 
   
As $M$ is a free $\O\otimes W$-module, we write $2g=\rank_W M$ and $g=fr$.
\begin{lemma}\label{42}
  Any quasi-polarized \dieu $\O$-module $\ul M$ has $2r$ slopes
  $\lambda_1\le\lambda_2\le \dots \le \lambda_{2r}$ with each slope
  multiplicity $f$. In other words, the Newton polygon of $M$
  satisfies the usual symmetric and integral conditions, and an 
  additional condition that the $x$-coordinate of the breaking points
  $(x,y)$ is divisible by $f$. 
\end{lemma}
\begin{proof}
  This should be well-known to experts; we supply a proof for
  convenience of the reader. Consider the iteration $F^f$ of the
  Frobenius operator on $M$, which preserves each factor $M^i$. 
  It follows from Katz's sharp slope estimate \cite{katz:slope} that
  the lowest slopes of $F^f$ on each $M^i$ are the same. Taking the
  exterior powers $\bigwedge^\bullet M^i$, each $M^i$ has the same slope
  sequence for $F^f$. Therefore $M$ has $2r$ slopes with each slope of
  multiplicity $f$.  
\end{proof}

\subsection{Example.}
\label{43}
When $f=3$ and $r=2$, the possible slope sequences of $M$
have the form  $\{\lambda_1^3,\lambda_2^3, \lambda_3^3, \lambda_4^3\}$
(with each multiplicity $3$) and are listed below:
\[ (\lambda_1,\lambda_2,\lambda_3,\lambda_4)=
 \left (\frac{1}{2},\frac{1}{2},\frac{1}{2},\frac{1}{2}\right )  
 > \left (\frac{1}{3},\frac{1}{2},\frac{1}{2},\frac{2}{3}\right )  
 > \left (\frac{1}{3},\frac{1}{3},\frac{2}{3},\frac{2}{3}\right ) 
 >\] \[ >  \left (0,\frac{1}{2},\frac{1}{2},1\right ),\  
   \left (\frac{1}{6},\frac{1}{6},\frac{5}{6},\frac{5}{6} \right )  
 > \left (0,\frac{1}{3},\frac{2}{3},1\right ) 
 > \left (0,0,1,1\right ). \]

\subsection{}
\label{44}

Recall the {\it $a$-number} of a \dieu module $M$, denoted by $a(M)$, 
is the dimension of the $k$-vector space $M/(F,V)M$.
We call $M$ \emph{local-local} if the operators $F$ and $V$ are 
nilpotent on $M/pM$. For a \dieu $\O$-module $\ul M$,  
the \emph{$a$-type} of $\ul M$ is defined as follows   
\[ \ul a(\ul M):=(a_i)_{\iz}, \quad a_i:=\dim_k (M/(F,V)M)^i. \]

\begin{thm} \label{T:normal}
Let $\ul M$ be a local-local separably quasi-polarized \dieu $\O$-module over
$k$ such that $\ul a(\ul M)=(1,0,\dots ,0)$.
Then there exists a symplectic
basis $X_0,X_1,\dots,X_{g-1},Y_0,Y_1,\dots,Y_{g-1}$ such that 
$X_i\in M^i$ and $Y_i\in(VM)^i$ for all $i=0,\dots ,{g-1}$, 
and the Frobenius operator $F$ has the following normal form:
\[ 
\begin{pmatrix}
  0 & 0 & \dots & 0 & & 0        & 0        & \dots & 0          & -p \\
  1 & 0 & \dots & 0 & & pa_{2,2} & pa_{2,3} & \dots & pa_{2,g} &  0 \\
  \vdots&\ddots&\ddots&\vdots& &\vdots&\ddots&\ddots&\vdots&\vdots \\
  0 & \dots & 1 & 0 & & pa_{g,2} & pa_{g,3} & \dots & pa_{g,g} & 0 \\
    &       &   &   & &          &          &       &    & \\     
  0 & \dots & 0 & 1 & & 0 & 0 & \dots &  \dots & 0 \\
  0 & \dots &\dots& 0 & & p& 0 &\dots &\dots &0 \\
  0 & \ddots&\ddots& 0& & 0& p &\dots &\dots &0 \\
  \vdots&\ddots&\ddots&\vdots& &\vdots&\vdots&\ddots&\vdots&\vdots\\
  0 & \dots & \dots & 0 & & 0 & \dots & \dots & p & 0\\ 
\end{pmatrix}
\]  
where $a_{i,j}=0$  if  $i\not\equiv j$ (mod $f$) and
$a_{i,j}=a_{j ,i}$ for all $i,j$ between $2$ and $g$. 
\end{thm}

\begin{proof}
We prove the following statement by induction 
\begin{itemize}
\item [($I_n$):] There is an $X\in M^0$ and $Y_i\in (VM)^i$ for
  $i=0,\dots, g-1$ such that 
\[ Y_0 \equiv F^g X\quad (\text{mod } p^n M),\quad  
Y_{g-1} \equiv -VX\quad (\text{mod } p^n M), \]
and if we set $X_i:=F^iX$ for $i=0,\dots,g-1$, then $\{X_i,Y_i\}$ is a
symplectic basis modulo $p^n M$. 
\end{itemize}

As $a(M)=1$, $M/pM$ is generated by one element $x$ over $k[F,V]$; it
is generated over $k$ by $x, F^1 x, \dots, F^{g-1} x, Vx,\dots,
V^{g-1}x$, and $F^g x$ (or $V^gx$). Particularly $F^g(M/pM)$ is
one-dimensional. 
Since the polarization is separable, it induces a perfect pairing 
\[ M/(F,V)M\times F^g(M/pM)\to W/pW\simeq k. \] 
Thus there exists an element $X\in M^0$ such that $\<X,F^g
X\>\not\equiv 0$ mod 
$p$. As $k$ is algebraically closed, we may adjust 
$X$ by a scalar multiple such that $\<X,F^g X\>\equiv 1$ mod
$p$. This yields $\<X,-V^gX\>\equiv 1$ mod $p$. 

As $-V^g\ol {X}\in F^g\ol {M}$, where $\ol{M}:=M/pM$, we obtain $-V^g
X\equiv F^g X$ mod $pM$. Define
\begin{equation} \label{E:adjust}
 Y_0:=F^g X,\quad Y_i:=-V^{g-i}X, \quad X_i:=F^i X \quad
 \forall\, i=1,\dots, g-1.
\end{equation}
Then we have $X_i\in M^i, Y_i\in (VM)^i$ and 

\begin{equation}
  \label{eq:342}
  \<X_i,Y_j\>\in
   \begin{cases}
   pW    &\text{if $i>j$,} \\
   1+pW  &\text{if $i=j$,} \\
   W     &\text{otherwise.}
   \end{cases}
\end{equation}

Set $X'=X+a_1 Y_f+\dots +a_{r-1}Y_{(r-1)f}$. It follows from
(\ref{eq:342}) that the $(r-1)\times (r-1)$ matrix
$(\<Y_{if},X_{jf}\>)_{1\le i, j<r}$ is 
non-degenerate mod $p$; we can find solutions $a_i\in W$ for
$1\le i <r$ such that $\<X',X_{jf}\>\equiv 0$ mod $pW$ for 
$j=1,\dots, r-1$. 
Note that 
$F^iX'\equiv F^iX$ mod $pM$. Replacing $X$ by $X'$, we have $\<X,F^i
X\>\equiv 0$ mod $pW$ for $0\le i<g$. 
Set $X_0:=X$ and re-define $Y_i$ and $X_i$ by \eqref{E:adjust}. Thus
we have verified $n=1$ of the statement $(I_n)$.

  
Suppose that we have a basis $\{X_i,Y_i\}$ satisfying $(I_n)$, we want
to find a new basis $\{X'_i,Y'_i\}$ that satisfies $(I_{n+1})$ and
that $X'_i\equiv X_i,\, Y'_i\equiv Y_i$ (mod $p^n M$). Let 
\[ X':=X+a_1 p^n Y_f+a_{2}p^n Y_{2f}+\dots+a_{(r-1)}p^n Y_{(r-1)f}.
\]
Since $F^i X'\equiv F^i X$ mod $p^{n+1}M$ for $i>0$, we have $\<X',F^i
X'\>\equiv \<X',X_i\>$ mod $p^{n+1} W$. Dividing $\<\, ,\, \>$ by
$p^n$, we get a system of linear equations mod $p$ with coefficients
$\<Y_j,X_i\>$. Therefore we can  find $a_i$, unique modulo $p$, such
that $\<X',F^i X'\>\equiv 0$ mod $p^{n+1}W$ for $i=1,\dots, g-1$. 
Clearly, $X'\equiv X$ mod $p^n M$. 

Suppose $\<X',F^g X'\>\equiv 1+a p^n$ mod $p^{n+1}W$. Choose $b\in W$
so that $b+b^{\sigma^g}+a\equiv 0$ mod $p$. Set $X'':=(1+b p^n) X'$, then
we have $\<X'', F^g X''\>\equiv 1$ mod $p^{n+1}W$ from the choice of
$b$. Replacing $X'$ by
$X''$, we find an $X'$ satisfying 
\begin{equation*}
  \begin{split}
  \<X', F^iX'\> & \equiv 0\quad \text{mod } p^{n+1}W\quad
  \text {for } 0<i<g, \\
  \<X',F^g X'\> & \equiv 1\quad \text{mod } p^{n+1}W, 
  \end{split}
\end{equation*}
and $X'\equiv X$ mod $p^nM$.

Set $X'_i:=F^i X'\!,\ Y'_{g-1}:=-VX'\!,\ Y'_0:=F^g X'$.  We can easily
adjust the rest of $Y_i$ by adding some linear combinations of
$p^nY_i$, and get $Y'_i$ for $i=1,\dots, g-2$ such that
$\{X'_i,Y'_i\}$ satisfies $(I_{n+1})$. 

By completeness of $M$, we get a symplectic basis $\{X_i,Y_i\}$ with
\[ (*):\quad  X_i\in M^i,\,\, Y_i\in (VM)^i,\,\, X_i=F^iX_0, \,\, Y_0=F^g X_0,
\text{ and } Y_{g-1}=-VX. \]

It follows from $(*)$ that $FY_j=pY_{j+1}+\sum_{i=1}^{g-1} p
a_{i+1,j+2}X_i$ for $j=0,\dots, g-2$ such that $a_{i+1,j+2}=0$ except
$i+1\not\equiv j+2$ mod $f$. 
That is,  $a_{i,j}\not= 0$ only when $i\equiv j$ (mod $f$). \qed

\end{proof} 

\section{Local moduli of the Newton strata}
\label{sec:05}

\subsection{}
\label{51}
We review the method of Cayley-Hamilton developed in Oort
\cite{oort:cayley}. Let $M$ be a \dieu module of rank $2g$. Suppose 
that the representative matrix of the Frobenius operator on $M$ with
respect to a basis has the following form (cf. Theorem~\ref{T:normal})
\[ \begin{pmatrix}  
B&C&E\\
0&D&0 \end{pmatrix}, \]
where $B$ is the $g\times (g-1)$ matrix with entries
$b_{i+1,j}=\delta_{ij}$, $D$ is the $g\times g$ diagonal matrix $\mathrm{diag}
(1,p,p,\dots,p)$, $E$ is the column vector $(-p,0,\dots,0)^t$ 
and $C$ is a $g\times g$ matrix of the form
\begin{equation}\label{eq:511}
   \begin{pmatrix}
   a_{1,1}&pa_{1,2}&\dots &pa_{1,g}\\
   a_{2,1}&pa_{2,2}&\dots &pa_{2,g}\\
   \vdots & \vdots &\ddots&\vdots \\
   a_{g,1}&pa_{g,2}&\dots&pa_{g,g}
   \end{pmatrix}. 
\end{equation}
The matrix $C$ is the main part of the representative matrix.
The assertion of this method says that the Newton polygon of $M$ can
be read off from that of the following polynomial 
\cite[Proposition 2.7]{oort:cayley} 
\[ X^{2g}+\sum_{k=g-1}^{0}\sum_{i=1}^{g-k}p^{i-1}a_{i+k,i}^{\sigma^{g-i+1}} 
   X^{g+k} +\sum_{k=1}^{g-1}\sum_{i=1}^{g-k} 
   p^{i+k-1}a_{i,i+k}^{\sigma^{g-i-k+1}} X^{g-k}-p^g. \] 
Suppose that every $a_{i,j}$ is either a unit, or sufficiently close
  to $0$ in $W$. In this case the valuation of the coefficient of
   every $X^{g+k}$ (resp. every $X^{g-k}$)
   is $\ord(p^{i-1}a_{i+k,k})$ (resp. $\ord(p^{i+k-1}a_{i,i+k})$) 
   for the smallest $i$ such that
   $a_{i+k,i}$ (resp. $a_{i,i+k}$) is a unit. 
  Then the Newton polygon of $M$ is the lower convex
  polygon obtained by removing from the following digram the
   $a_{i,j}$ with $p|a_{i,j}$ 

\begin{equation}\label{eq:512}
   \begin{matrix}
   &&&&&&&&&-1\\
   &&&&& a_{g,g}&\dots&\dots&a_{1,g}\\
   &&&&\cdot &\vdots&\cdot&\cdot&\\
   &&&\cdot&\cdot &\vdots&\cdot&\\
   &&1&a_{g,1}&\dots&a_{1,1} &&&\\
   \end{matrix}
\end{equation} 
where numbers are placed on the lattices of a $2g\times g$ box, $-1$
(resp. $1$) is placed at $(2g,g)$ (resp. $(0,0)$).

\subsection{}
\label{52}
\ncmd{\Cart}{\mathrm{Cart}_p}
We recall the deformation theory of \dieu modules developed in
Norman \cite{norman:algo} and Norman-Oort \cite{norman-oort}. 
We follow the convenient setting of \cite[Section 0]{norman:algo} and 
\cite[Section 2, p.~1011]{chai-norman:gamma2}.

Let $R$ be a commutative ring of \ch $p$. Let $W(R)$ denote the ring of
Witt vectors over $R$, equipped with the Verschiebung $\tau$ and
Frobenius $\sigma$:
\begin{equation*}
  \begin{split}
    (a_0,a_1,\dots)^\tau&=(0, a_0, a_1,\dots) \\
    (a_0,a_1,\dots)^\sigma&=(a_0^p, a_1^p,\dots).
  \end{split}
\end{equation*}
Let $\Cart(R)$ denote the Cartier ring $W(R)[F][\![V]\!]$ modulo the
relations
\begin{itemize}
\item $FV=p$ and $VaF=a^\tau$,
\item $Fa=a^\sigma F$ and $Va^\tau=aV, \ \forall\, a\in W(R).$
\end{itemize}

A left $\Cart(R)$-module is \emph{uniform} if it is complete and
  separated in the $V$-adic topology. A uniform $\Cart(R)$-module $M$ is
  \emph{reduced} if $V$ is injective on $M$ and $M/VM$ is a free
  $R$-module. A \emph{\dieu module over $R$} is a finitely generated
  reduced uniform $\Cart(R)$-module.

There is an equivalence of categories between the category of finite
dimensional commutative formal groups over $R$ and the category of
\dieu modules over $R$. We denote this functor by $\bfD_*$. The tangent
space of a formal group $G$ is canonically isomorphic to
$\bfD_*(G)/V\bfD_*(G)$.

\subsection{}
\label{53}
\def\k{{\rm k}}
Let $\ul M_0=(M_0,\<\, ,\, \>,\iota)$ be a local-local separably
quasi-polarized \dieu  
$\O$-module with $a(M_0)=1$. For simplicity we assume
that $\ul a(\ul M_0)=(1,0,\dots ,0)$. Let $\{X_i,Y_i\}$ be a 
symplectic basis as in Theorem~\ref{T:normal}. Set
\begin{equation}\label{eq:53}
  \begin{split}
    e^\ell_\k & :=X_{(\k-1)f+\ell}\quad \text{ for \ } 0\le \ell\le f-1, 
    1\le \k \le r, \\
    f^\ell_\k & :=Y_{(\k-1)f+\ell}\quad \text{ \ for } -1\le \ell\le f-2, 
    1\le \k \le r, \\
  \end{split}
\end{equation}
with the notation $X_{i\pm g}=X_{i}$ and $Y_{i\pm g}=Y_i$. 
Then $\{e^i_\k, f^i_\k\}_{1\le \k \le r}$ is a symplectic
  $W$-basis of $M^i$ for each $\iz$. It is convenient to extend the
  notation $f^\ell_\k$ and $e^\ell_\k$ for all integers $\ell,\k$ by the rule
  (\ref{eq:53}).  
We have 
\begin{align*}
 Fe^\ell_j\quad & =e^{\ell+1}_j\quad\quad \text{except $(\ell,j)=(f-1,r)$,} \\
 Fe^{f-1}_r&=f^0_1, \\
 Ff^\ell_j\quad &=pf^{\ell+1}_j+\sum_{i=1}^rpa^{\ell+1}_{i,j}e^{\ell+1}_i
 \quad\ \text{except $(\ell,j)=(-1,1)$,}\\
  Ff^{-1}_{1} &=-pe^0_1. 
\end{align*}


\subsection{}
\label{54}

\newcommand{\Def}{\mathrm{Def}}

Let $\Def[\ul M_0]$ be the deformation functor of $\ul M_0$. 
This functor is smooth; moreover,  for any local Artinian $k$-algebra
$(R,\grm)$ with residue field $k$ there is a bijection \cite{norman:algo}
\[ \Def[\ul M_0](R)\simeq \Hom_\O (VM_0/pM_0,
M_0/VM_0)^{\rm sym}\otimes_k \grm. \]
Here we identify $VM_0/pM_0$ with the dual $(M_0/VM_0)^*$ by the
pairing $\<\, ,\,\>$. Let $\{(e^\ell_{\k})^*, (f^\ell_\k)^*\}$ be the dual
basis for the basis $\{e^\ell_\k, f^\ell_\k\}$. 
We may identify $(e^\ell_{\k})^*$ with $f^\ell_\k$ and $(f^\ell_\k)^*$ 
with $e^\ell_\k$. Thus any homomorphism that sends
$\bar{f}^\ell_j\mapsto \sum t^{\ell}_{i,j}\bar{e}^{\ell}_i$ with
$t^\ell_{i,j}=t^\ell_{j,i}\in{\mathfrak m}$ gives rise to 
a deformation of $\ul M_0$  over $R$; and conversely any 
deformation of $\ul M_0$ over $R$ arises from such a unique 
homomorphism. Hence
\[ R_{\rm univ}:=k[|t^\ell_{i,j}|]/(t^\ell_{i,j}-t^\ell_{j,i})= 
  k[|t^\ell_{i,j}|]_{0\le \ell <f,\ 1\le j\le i\le r}\] 
is the universal deformation ring of $\ul M_0$.

By \cite{norman:algo} or \cite{norman-oort}, 
the universal deformation $M_{R_{\rm univ}}$ of $\ul M_0$ over
$R_{\rm univ}$ is given by the following generators and relations:
\begin{equation}\label{eq:54}
\begin{split}
 F{e}^\ell_j\ \ &={e}^{\ell+1}_j\quad\text{except
  $(\ell,j)=(f-1,r)$,}  \\
 F{e}^{f-1}_r&={f}^0_1+\sum_{i=1}^{r}T^0_{i,1}
 {e}^0_i, \\
 {f}^\ell_j\ \ &=V\left({f}^{\ell+1}_j+\sum_{i=1}^r
 (a^{\ell+1}_{i,j}+ T^{\ell+1}_{i,j}){e}^{\ell+1}_i \right)
 \quad\text{except $(\ell,j)=(-1,1)$,}\\
 {f}^{-1}_{1}&=-V{e}^0_1, 
\end{split}
\end{equation}
where 
$T^\ell_{i,j}:=(t^\ell_{i,j},0,\dots)$ is the 
Teichm\"{u}ller lifting of $t^\ell_{i,j}$ for each $\ell,i,j$.

\subsection{}
\label{55}
  Let $\ul G_0$ be the $p$-divisible group with 
  additional structures corresponding to $\ul M_0$ and $\ul G^{\rm
  univ}$ be its universal deformation over $\calU:=\Spec\, R_{\rm univ}$. 
  Let $\beta_0=NP(M_0)$ and $\beta$ be an admissible Newton
  polygon below $\beta_0$. Denote by $\calU^{\ge \beta}$ the reduced
  closed subscheme of $\calU$ whose points has Newton polygon lying
  over or equal to $\beta$. 
  By the theory of Cayley-Hamilton, the closed Newton stratum $\calU^{\ge
  \beta}$ is defined by the equations $t^\ell_{i,j}$ whose
  corresponding ``position'' in the Oort diagram ((\ref{eq:512}) of
  (\ref{51})) is strictly below the Newton polygon $\beta$; see the
  discussion in Oort \cite[Section 3]{oort:cayley}; also see an
  example in (\ref{57}). 

  From (\ref{eq:511}) of (\ref{51}) and (\ref{eq:54}) of (\ref{54}), it
  is not hard to see that the corresponding position of
  $t^\ell_{\k+j,j}$ is $\left( (r-\k)f, (j-1)f+\ell \right )$. Put
\begin{equation}\label{eq:55}
 S(\beta):=\{\,t^\ell_{i,j}\ ;\ j\leq i \text{ and
the position of $t^\ell_{i,j}$ is higher or on $\beta$}\ \},
\end{equation}
thus $\calU^{\ge \beta} \simeq \Spec \,k[|t|]_{t\in S(\beta)}$.

When $\beta_0$ is supersingular, $\calU^{\ge \beta_0}$ is the
supersingular locus in $\calU$. One computes that 
\[ \dim \calU^{ss} = |S(\beta_0)|=
\sum_{i=1}^r[i\,f/2]=
  \begin{cases}
  \k\frac{r(r+1)}{2} &\text{if $f=2\k$,}\\
  \k\frac{r(r+1)}{2}+[r^2/4] &\text{if $f=2\k+1$.}
  \end{cases} \]

\begin{thm}\label{56} 
Let $\ul M_0$ be a separably quasi-polarized \dieu $\O$-module with
$a(M_0)=1$ and $\beta\le NP(M_0)$ be an admissible Newton polygon. 
Then the closed Newton stratum $\calU^{\ge \beta}$ of the universal 
deformation $\calU$ is formally smooth of dimension $|S(\beta)|$
above; the generic point has Newton polygon equal to $\beta$. 
If $M_0$ is supersingular, the supersingular locus $\calU^{ss}$ in
$\calU$ has dimension 
\[ \dim \calU^{ss}=
  \begin{cases}
  \k\frac{r(r+1)}{2} &\text{if $f=2\k$,}\\
  \k\frac{r(r+1)}{2}+[r^2/4] &\text{if $f=2\k+1$.}
  \end{cases} \]
\end{thm}

  We have done the case where $M_0$ is local-local.
  When $M_0$ is not local-local, we have a decomposition 
  $M_0=M_0^\text{\,loc-loc}\oplus M_0^{\rm ord}$. 
  Write $\calU_1$ for the universal 
  deformation for $M_0^\text{\, loc-loc}$ and $\beta=NP(M_0^{\rm
  ord})\cup \beta_1$. The defining equations for
  $\calU^{\ge \beta}$ are the same as those for $\calU_1^{\ge
  \beta_1}$. Therefore, we obtain the same dimension formula as
  before. 

\subsection{Example.}
\label{57}
When $f=3,\ r=2$, the main part of the representative matrix of the
Frobenius operator on the universal \dieu module is
\[ 
\begin{pmatrix}
  T^0_{1,1} &&& p(a^0_{1,2}+T^0_{1,2})\\
    &p(a^1_{1,1}+T^1_{1,1})&&& p(a^1_{1,2}+T^1_{1,2}) \\
&&p(a^2_{1,1}+T^2_{1,1})&&& p(a^2_{1,2}+T^2_{1,2})\\
T^0_{2,1}&&& p(a^0_{2,2}+T^0_{2,2})\\
& p(a^1_{2,1}+T^1_{2,1})&&& p(a^1_{2,2}+T^1_{2,2})\\
&&p(a^2_{2,1}+T^2_{2,1})&&& p(a^2_{2,2}+T^2_{2,2})\\
\end{pmatrix} \]
We have the following picture 

\epsfig{file=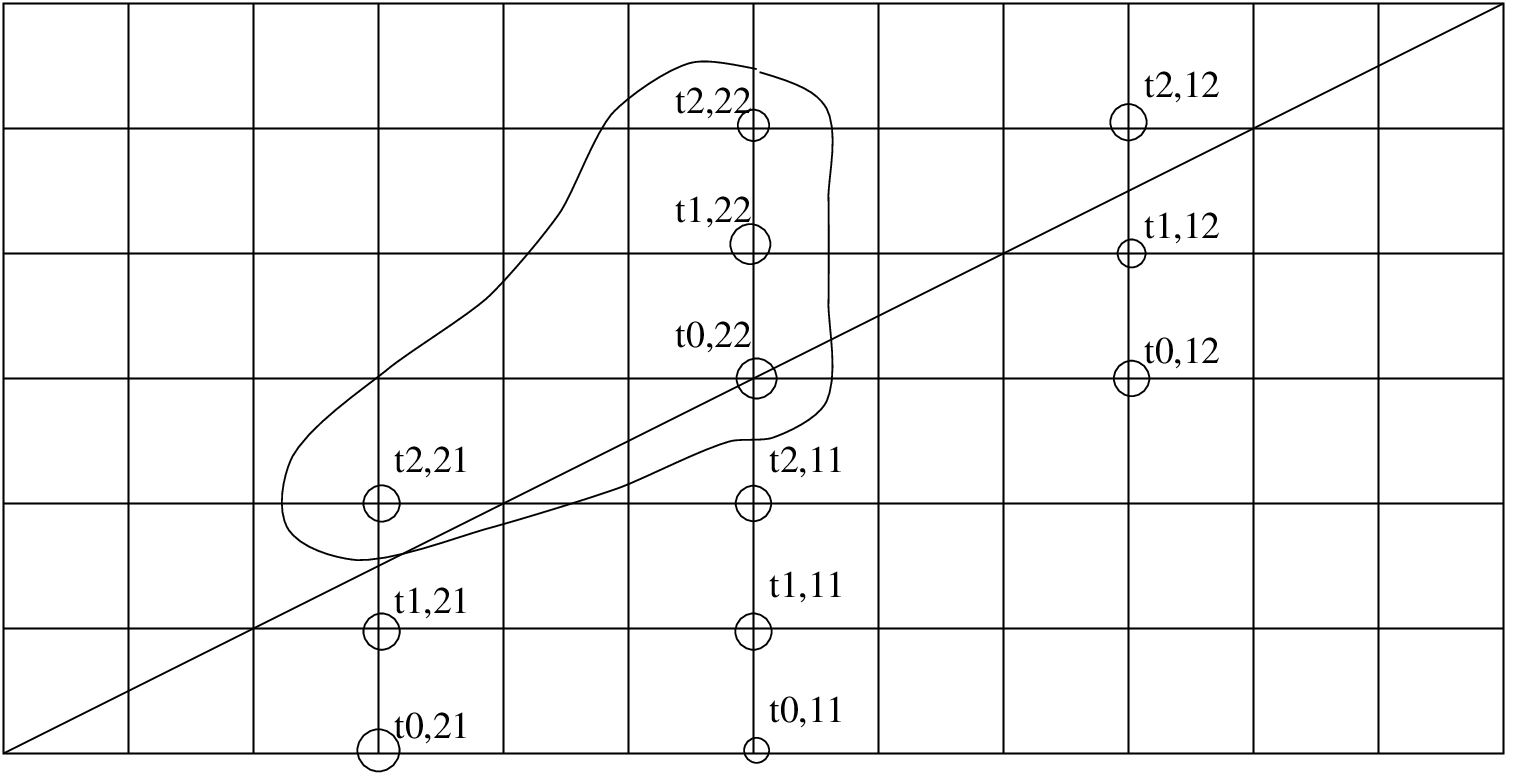,angle=270,width=4in}

\ \\ \\
Here $t2,11$ denotes $t^2_{1,1}$ etc.
Particularly $\calU^{ss} \simeq \Spec\, k[|t^2_{2,1},t^0_{2,2},
t^1_{2,2},t^2_{2,2}|]$ when $M_0$ is supersingular.

\subsection{}
\label{58}
We now return to complete Theorem~\ref{12}--\ref{14}. Using the
Serre-Tate theorem, the deformation problem is reduced to that for
$p$-divisible groups with the additional structure. As $p\nmid
\Delta$, one reduces to that for quasi-polarized $p$-divisible
$\O$-groups by the Morita equivalence. By Theorem~\ref{56},
Theorem~\ref{12} is proved. 

Using Theorem~\ref{11}, one has a separably polarized supersingular
point $\ul A_0$ in the open subset $U$. Then Theorem~\ref{13} follows
from Theorem~\ref{12} immediately by induction.

Theorem~\ref{14} follows from Theorem~\ref{13} and Theorem~\ref{11} (2).

\bigskip
\npr Institute of Mathematics \\ Academia Sinica \\
128 Academia Rd.~Sec.~2, Nankang \\ Taipei, Taiwan 115 \\
E-mail Address: chiafu@math.sinica.edu.tw

\end{document}